\documentclass[10pt,reqno]{amsart}
\usepackage{amssymb}
\usepackage[all]{xy}

\theoremstyle{plain}
\newtheorem{prop}{Proposition}
\newtheorem{theo}{Theorem}
\newtheorem{coro}{Corollary}
\newtheorem{lemm}{Lemma}
\theoremstyle{remark}
\newtheorem{rema}{Remark}
\theoremstyle{definition}

\newtheorem{exam}{Example}
\newtheorem{nota}{Notation}

\newcommand{\eqto}{\stackrel{\lower1.5pt\hbox{$\scriptstyle\sim\,$}}\to}

\newcommand{\al}{\alpha}

\def\lra{\longrightarrow}
\newcommand{\vlra}{\relbar\joinrel\longrightarrow}
\newcommand{\vvlra}{\relbar\joinrel\vlra}
\newcommand{\vvvlra}{\relbar\joinrel\vvlra}
\newcommand{\vvvvlra}{\relbar\joinrel\vvvlra}
\newcommand{\vvvvvlra}{\relbar\joinrel\vvvvlra}

\newcommand{\clap}[1]{\hbox to 0pt{\hss #1\hss}}

\def\cC{{\mathcal C}}

\def\cE{{\mathcal E}}

\def\cO{{\mathcal O}}

\def\Ker{{\rm Ker}}
\def\Gal{{\rm Gal}}

\def\ra{\rightarrow}

\def\C{{\mathbb C}}

\def\PP{{\mathbb P}}
\def\Q{{\mathbb Q}}
\def\R{{\mathbb R}}

\def\Z{{\mathbb Z}}
\def\C{{\mathbb C}}

\def\Pic{{\rm Pic}}

\def\Hom{{\rm Hom}}
\def\Br{{\rm Br}}

\begin{document}
\title[Del Pezzo surfaces of degree two]
{On the arithmetic of del Pezzo surfaces of degree~2}
\author{Andrew Kresch}
\address{
  Mathematics Institute,
  University of Warwick,
  Coventry CV4 7AL,
  United Kingdom
}
\email{kresch@maths.warwick.ac.uk}
\author{Yuri Tschinkel}
\address{
  Mathematisches Institut,
  Georg-August-Universit\"at G\"ottingen,
  Bunsenstrasse 3-5,
  D-37073 G\"ottingen, Germany
}
\email{yuri@uni-math.gwdg.de}

\dedicatory{To Sir Peter Swinnerton-Dyer}
\date{23 February 2004}
\thanks{The first author was supported by
an EPSRC Advanced Research Fellowship.
The second author was supported by the NSF grant 0100277.}

\begin{abstract}
We study the arithmetic of certain del Pezzo surfaces of degree 2. 
We produce examples of Brauer-Manin obstruction to the Hasse principle,
coming from $2$- and $4$-torsion elements in the Brauer group.
\end{abstract}
\maketitle
\tableofcontents

\setcounter{section}{0}
\section{Introduction}
\label{sect:introduction}
Del Pezzo surfaces are smooth projective surfaces, isomorphic over the
algebraic closure of the base field to $\PP^1\times\PP^1$ or
the blow-up of $\PP^2$ in up to 8 points in general position.
In the latter case the del Pezzo surface has degree equal to 9 minus
the number of points in the blow-up.
The arithmetic of del Pezzo surfaces over number fields is an active
area of investigation.
It is known that the Hasse principle holds for del Pezzo surfaces of
degree at least 5.

Counterexamples to the Hasse principle were discovered for
del Pezzo surfaces of degrees 3 and 4 (see \cite{swinnerton-dyer62}
and \cite{bsd75}, respectively).
A growing body of evidence (for instance, \cite{CCS}) led to the
question of whether the failure of the Hasse principle for del Pezzo surfaces
is always explained by the Brauer-Manin obstruction;
this question is specifically raised by
Colliot-Th\'el\`ene and Sansuc in \cite{angers}.
Computer verifications for diagonal cubics in \cite{cubdiag}
and theoretical advances, such as
\cite{stn89}, \cite{salberger}, \cite{swinnerton-dyer01},
lend support to an affirmative answer to this question.

A del Pezzo surface of degree 2
can be realised as a double cover of $\PP^2$ ramified 
in a smooth quartic curve. In this note we consider surfaces $S$
over $\Q$ of the form 
\begin{equation}
\label{eqn:d2}
w^2=Ax^4+By^4+Cz^4.
\end{equation}
We compute the Galois-theoretic invariant $\Br(S)/\Br(\Q)$ and
produce examples of obstruction to the Hasse principle
(see \cite{cubdiag}, \cite{manin} for background).
We obtain:
\begin{theo}
\label{th:main}
Let $S$ have the form \eqref{eqn:d2}, where $A$, $B$, and $C$ denote
nonzero integers.
Then $\Br(S)/\Br(\Q)$ is isomorphic to one of the following groups:
\begin{gather*}
(1), \ \Z/2, \ \Z/4, \ (\Z/2)\oplus(\Z/2), \\
(\Z/4)\oplus(\Z/2), \ (\Z/2)\oplus(\Z/2)\oplus(\Z/2).
\end{gather*}
\end{theo}
The simplest of our examples, the case $p=3$ of Example \ref{exam:obstr2},
is the assertion that
\begin{equation}
\label{eqn:-6-32}
w^2=-6x^4-3y^4+2z^4
\end{equation}
has no rational solutions aside from the trivial solution.
Here, a completely down-to-earth formulation of the proof is that,
rewriting \eqref{eqn:-6-32} as
$$w^2+(2x^2-y^2)^2=2(x^2+y^2+z^2)(-x^2-y^2+z^2)$$
and supposing $(w,x,y,z)$ to be an integer solution with no common
prime factors, we get a contradiction to the expression on the right being
a sum of squares from the factors $x^2+y^2+z^2$ and $-x^2-y^2+z^2$
each being congruent to
$3$ mod $4$, yet having no common prime factor congruent to $3$ mod $4$.
Example \ref{exam:obstr2} shows that \eqref{eqn:-6-32} fits into an
infinite sequence of counterexamples to the Hasse principle.
A more sophisticated example, Example \ref{exam:obstr5}, is of
particular interest, since
the obstruction comes from a $4$-torsion element in the Brauer group.
By \cite{swinnerton-dyer93}, only $2$- and $3$-torsion Brauer group elements
occur for del Pezzo surfaces of degree $\ge 3$.

The tool we use is group cohomology.
Let $F$ be a Galois extension of $\Q$, and let $G$ denote the
Galois group $\Gal(F/\Q)$.
If $\Pic(S_F)$ is equal to the geometric Picard group
$M:=\Pic(S_{\overline\Q})$ then we have
\begin{equation}
\label{eqn:H1}
\Br(S)/\Br(\Q)=H^1(G,M).
\end{equation}
More generally, the Hochschild-Serre spectral sequence gives rise to the
following exact sequence:
\begin{multline}
\label{eqn:HS}
0\lra\Pic(S)\lra \Pic(S_F)^G\lra\ker(\Br(\Q){\to}\Br(F))\\
\lra\ker(\Br(S){\to}\Br(S_F))\lra H^1(G,\Pic(S_F))\lra H^3(G,F^*).
\end{multline}
In this paper, we compute the group \eqref{eqn:H1} and
represent lifts of elements to $\Br(S)$ by
Azumaya algebras.
By \eqref{eqn:HS} and cohomological dimension,
such lifts exist after perhaps enlarging $F$;
what happens in practice is that it is often possible to take
$[F:\Q]$ quite small and still have
$H^1(G,\Pic(S_F))$ isomorphic to $\Br(S)/\Br(\Q)$ and
the final map in \eqref{eqn:HS} trivial.
Lastly, we explain the computation of local invariants and obtain
the above-mentioned examples.

In an Appendix we show that in the case of the
diagonal cubic surfaces considered in \cite{cubdiag} the present techniques
give rise to cyclic Azumaya algebras.
This simplifies the construction of cocycle representatives
and the local obstruction analysis, as compared with
the original consideration of bicyclic group cohomology.

We take a moment to highlight instances where the arithmetic of
del Pezzo surfaces of degree $2$ has already been studied.
Our examples are new, and this is the first systematic study of a
class of degree $2$ del Pezzo surfaces.
However, some special classes of degree $2$ del Pezzo surfaces do fall within
the scope of existing results.
These include:
\begin{itemize}
\item[(i)] \emph{Blow-ups of higher degree del Pezzo surfaces.}
One can, for instance, start with a degree $4$ del Pezzo surface which
violates the Hasse principle, such as can be found in \cite{bsd75},
and blow up a conjugate pair of points
(which, as mentioned in \cite{colliot-thelene84},
always exist on such a surface) to obtain a del Pezzo surface of
degree $2$ which is a counterexample to the Hasse principle.
\item[(ii)] \emph{Double coverings of Ch\^atelet surfaces.}
We are grateful to Colliot-Th\'el\`ene for providing the following example.
We consider the equation
\begin{equation}
\label{eq.chatelet}
w^2+x^2=(y^2-2)(3-y^2).
\end{equation}
This defines a (generalised) Ch\^atelet surface which fails to
satisfy the Hasse principle.
We now replace $x^2$ by $x^4$ in the equation;
there remain points in all completions of $\Q$, so we get
a degree $2$ del Pezzo surface which fails to satisfy the Hasse principle.
Note (cf.\ \cite{CCS}) that \eqref{eq.chatelet}
belongs to an infinite sequence of
counterexamples to the Hasse principle given by Iskovskih \cite{isk}.
\item[(iii)] \emph{Birational models of conic bundles with six
degenerate fibres.}
Many del Pezzo surfaces of degree $2$ fit this description;
the referee is credited with suggesting this source of examples.
Notably, for surfaces of the form
\begin{equation}
\label{eq.conicsix}
r^2+s^2 = f_2(t)f_4(t)
\end{equation}
with $f_2(t)$ and $f_4(t)$ irreducible polynomials of degrees $2$ and $4$,
respectively, Swinnerton-Dyer has shown that the Brauer-Manin obstruction
is the only obstruction to the Hasse principle \cite{swinnerton-dyer99}.
In fact the same is true for weak approximation; see \cite{skorobogatov}.
This means that on any smooth projective model,
the rational points are dense in the set of adelic points not
obstructed by Brauer classes.
As a concrete example, the Brauer-Manin obstruction is
the only obstruction to weak approximation
for the del Pezzo surface given by
$$w^2=5x^2z^2-4y^4+26y^2z^2-30z^4-4x^3z-16xy^2z+24xz^3.$$
This is birational to the surface \eqref{eq.conicsix} with
$f_2(t)=t^2+3$ and $f_4(t)=t^4+t^2+2$:
$$r=2\frac{y}{z}+t^3,\qquad s=\frac{x}{z}-t\frac{y}{z}+2t^2,\qquad
t=\frac{x^2+4y^2-6z^2}{w+xy}.$$
\end{itemize}
Our examples are not covered by cases (i)--(iii).
We discuss this briefly at the end of Section \ref{sec:examples}.

The authors would like to thank J.-L. Colliot-Th\'el\`ene for helpful
discussions and correspondence.

\section{Geometry}
\label{sect:geom}
Consider the surface $S$ given by the equation
$$
w^2=Ax^4+By^4+Cz^4
$$
in the weighted projective space $\PP(2,1,1,1)$,
where $A$, $B$, and $C$ are nonzero integers.
It is a double cover of $\PP^2$, branched over the twisted Fermat quartic
curve
$$0=Ax^4+By^4+Cz^4.$$
Let $a,b,c$ denote some chosen 4-th roots of $A,B,C$, respectively.  
The $56$ exceptional curves on $S$ are the pre-images of the
bitangents to the quartic. 
These are given by the following equations
\begin{gather}
\delta ax + by =  0 , \,\,\,\,\,  
\delta by+cz = 0  ,  \,\,\,\,\,
\delta cz +ax = 0  , \,\,\,\,  \text{ where } \delta^4=-1 , \label{eqn:p} \\
\al ax+ \beta by + \gamma cz =0 \,\,\,\, (\al^4=\beta^4=\gamma^4=1) 
\label{eqn:pp} . 
\end{gather}
Multiplying the equation \eqref{eqn:pp} by a scalar doesn't change the
line it defines, so it is natural to index the line by an element 
$(\al,\beta,\gamma)\in \mu_4^3/\mu_4$.
Each bitangent lifts to a pair of exceptional curves in $S$:
for example, the pre-image of the line given by $\delta ax +by=0$
is the pair of curves with equations 
$$
w=\pm c^2z^2 \ .
$$ 
These will be denoted by $L_{z,\delta,\pm}$. 
There are $24$ exceptional curves lying over the lines in \eqref{eqn:p}. 
The pre-images of the lines in \eqref{eqn:pp} are given by
\begin{equation}
\label{eqn:wpm}
w=\pm \sqrt{2}(\al \beta abxy+\beta\gamma bcyz +\al\gamma acxz) \ . 
\end{equation}
The ambiguity $\pm $ is resolved by scaling the tuple $(\al,\beta,\gamma)$;
we denote by $L_{\alpha,\beta,\gamma}$ the pre-image
\eqref{eqn:wpm} with the sign taken to be $+$, so now
$(\al,\beta,\gamma)$ is considered to be in $\mu_4^3/\mu_2$.
We thus have the following description of the exceptional curves on $S$.
\begin{prop}
\label{prop:exc}
The $56$ exceptional curves on the del Pezzo surface \eqref{eqn:d2} are
as follows, where $a$, $b$, and $c$ denote chosen $4$-th roots of $A,B,C$:
$$
\begin{array}{ccl}
L_{z,\delta,\pm}\,:  &  \delta a x + by =0,  & w=\pm c^2z^2,
                                        \qquad (\delta^4=-1)\ , \\
L_{x,\delta,\pm}\,:  &  \delta by + cz =0,   & w=\pm a^2x^2,
                                        \qquad (\delta^4=-1)\ , \\
L_{y,\delta,\pm}\,:  &  \delta cz + ax =0,   & w=\pm b^2y^2,
                                        \qquad (\delta^4=-1)\ ,  \\ 
L_{\al,\beta,\gamma} \,: & \al ax + \beta by + \gamma cz =0 & 
w=\sqrt{2}(\al \beta ab xy + \beta \gamma bc yz + \al \gamma acxz ), \\
&&\qquad\qquad\qquad\qquad (\alpha,\beta,\gamma)\in \mu_4^3/\mu_2\ . 
\end{array}
$$
\end{prop}

Geometrically, the Picard group of $S$ has rank $8$.
We choose the basis indicated in the following statement.
\begin{prop}
\label{prop:PicS}
Let $S$ be the del Pezzo surface \eqref{eqn:d2}, and
set $\zeta=e^{\pi i/4}$.
Then the geometric Picard group $\Pic(S_{\overline\Q})$ is the free abelian
group on the generators
$$
\begin{array}{llll}
v_1 = [L_{x,\zeta,+}] &  v_2 = [L_{x,\zeta^3,-}]  
&  v_3   = [L_{y,\zeta,+}] & v_4 = [L_{y,\zeta^3,-}]   \\ 
v_5 = [L_{z,\zeta,+}]  &  v_6 = [L_{z,\zeta^3,-}]
&  v_7   = [L_{i,i,i}]  &
v_8 = [L_{z,\zeta^7,-}]+[L_{z,\zeta^3,-}]+[L_{i,i,i}]  \ .
\end{array}
$$
The class $v_i$ has self-intersection $-1$ for $i\le 7$ and
self-intersection $1$ for $i=8$.
The intersection number of $v_i$ and $v_j$ is $0$ for $i\ne j$.
The anticanonical class is
\begin{equation}
\label{eqn:anticanonical}
-K_S=-v_1-v_2-v_3-v_4-v_5-v_6-v_7+3v_8.
\end{equation}
The identities displayed in Table \ref{tab:cl} hold
in $\Pic(S_{\overline\Q})$; these, coupled with \eqref{eqn:anticanonical}
determine the class of any exceptional curve.
\end{prop}

\begin{proof}
Each exceptional curve
has self-intersection $-1$. 
Each pair of curves lying above a bitangent to the Fermat quartic
has intersection number 2.
Other intersection numbers are 0 or 1 and are readily determined.
In particular, the intersection numbers among the $v_i$ are as claimed,
and the $v_i$ span $\Pic(S_{\overline\Q})$.
The anticanonical class is the class of any pair of curves lying
above a bitangent to the Fermat quartic.
The anticanonical class and the classes listed in Table \ref{tab:cl}
are determined by computing intersection numbers with the $v_i$.
\end{proof}

\begin{table}
$$
\begin{array}{ccc}
{[L_{x,\zeta^5,+}]=-v_1-v_7+v_8} && {[L_{x,\zeta^7,-}]=-v_2-v_7+v_8} \\
{[L_{y,\zeta^5,+}]=-v_3-v_7+v_8} && {[L_{y,\zeta^7,-}]=-v_4-v_7+v_8} \\
{[L_{z,\zeta^5,+}]=-v_5-v_7+v_8} && {[L_{1,1,i}]=-v_2-v_3+v_8} \\
{[L_{1,1,-1}]=-v_5-v_6+v_8} && {[L_{1,1,-i}]=-v_1-v_4+v_8} \\
{[L_{1,i,1}]=-v_1-v_6+v_8} && {[L_{1,i,-i}]=-v_3-v_5+v_8} \\
{[L_{1,-1,1}]=-v_3-v_4+v_8} && {[L_{1,-1,-1}]=-v_1-v_2+v_8} \\
{[L_{1,-i,1}]=-v_2-v_5+v_8} && {[L_{1,-i,i}]=-v_4-v_6+v_8} \\
{[L_{i,1,1}]=-v_4-v_5+v_8} && {[L_{i,1,-i}]=-v_2-v_6+v_8} \\
{[L_{i,-1,-1}]=-v_3-v_6+v_8} && {[L_{i,-1,-i}]=-v_1-v_5+v_8} \\
{[L_{i,-i,1}]=-v_1-v_3+v_8} && {[L_{i,-i,-1}]=-v_2-v_4+v_8}
\end{array}
$$
\caption{Classes of the exceptional curves}
\label{tab:cl}
\end{table}

\section{Galois group - generic case}
\label{sect:gal-gr}
Let $G$ be the Galois group of the extension 
\begin{equation}
\label{eqn:F}
F:=\Q(\zeta,a^2,b/a,c/a)
\end{equation}
over $\Q$ (where $\zeta=e^{\pi i/4}$). 
The subextension 
$\Q(\zeta)/\Q$ corresponds
to a normal subgroup $H$ of index 4. 
The quotient group is the Klein four-group. 
In the {\em generic} case, we have
$|G|=128$.
The Galois group can be described as follows.
\begin{prop}
\label{prop:ga}
Let $A$, $B$, and $C$ be nonzero integers, with chosen
$4$-th roots $a$, $b$, and $c$, respectively.
Let $F$ be as in \eqref{eqn:F}, and suppose
the degree of $F$ over $\Q$ is $128$.
Then the Galois group $G_0=\Gal(F/\Q)$ is generated by elements
$$
\sigma, \tau, \iota_a,\iota_b,\iota_c
$$
which act by
$$
\begin{array}{|ccc|ccc|ccc|ccc|ccc|ccc|} 
\hline
&&&&  \sigma &&& \tau &&& \iota_a   &&& \iota_b    &&& \iota_c  & \\
\hline &
a^2 &&& a^2 
&&&   a^2 
&&&   -a^2 &&& a^2   &&& a^2 & \\
& b/a &&& b/a &&& b/a 
&&&   -ib/a 
&&&   ib/a &&& b/a & \\
& c/a &&& c/a &&& c/a 
&&&   -ic/a 
&&&   c/a &&& ic/a & \\
& \zeta &&& \zeta^{-1}  &&&
\zeta^3 &&&  
\zeta &&& \zeta &&& \zeta & \\ 
\hline
\end{array}
$$
The action on the exceptional curves is as follows:
$$
\begin{array}{|l|c|c|c|c|c|} 
\hline
 &  \sigma & \tau & \iota_a   & \iota_b    & \iota_c   \\
\hline
L_{z,\delta,s} & L_{z,\sigma(\delta),s} 
&   L_{z,\tau(\delta),s} 
&   L_{z,i\delta,s} & L_{z,-i\delta,s}   & L_{z,\delta, -s} \\
L_{x,\delta,s} & L_{x,\sigma(\delta),s} & L_{x,\tau(\delta),s} 
&   L_{x,\delta,-s} 
&   L_{x,i\delta,s} & L_{x,-i\delta,s}  \\
L_{y,\delta,s} & L_{y,\sigma(\delta),s} & L_{y,\tau(\delta),s} 
&   L_{y,-i\delta,s} 
&   L_{y,\delta,-s} & L_{y,i\delta,s}  \\
L_{\al,\beta,\gamma} & L_{\al^{-1},\beta^{-1}, \gamma^{-1}}  &
L_{i\al^{-1}, i\beta^{-1},i\gamma^{-1}} &  
L_{i\al,\beta,\gamma} & L_{\al ,i\beta,\gamma} & L_{\al,\beta,i\gamma}\\ 
\hline
\end{array}
$$
\end{prop}

\section{Group cohomology}
\label{sect:gr}
We start with a review.
If $G$ is a group, a standard free resolution of $\Z$ is
\begin{equation}
\label{eqn:1}
\cC_{\bullet}^G:= \cdots 
\Z[G\times G\times G]  \ra  \Z[G\times G]\ra \Z[G] \ , 
\end{equation}
where
the augmentation map $\Z[G]\ra \Z$ is given 
by $g\mapsto 1$ (for all $g\in G$) and 
where each map in $\cC_{\bullet}^G$ is of the form 
$$
(g_0,\ldots g_n)\mapsto
\sum_{i=0}^n (-1)^i (g_0,\ldots, \widehat{g_i},
\ldots, g_n) \ .
$$
The action of  $g\in G$ on any of the terms in \eqref{eqn:1}
is the diagonal left multiplication action.
We may identify 
\begin{equation}
\label{eqn:4}
\begin{array}{ccc}
\Z[G\times G] &\simeq    &  \bigoplus_{g'\in G}  \Z[G] \ ,  \\
(g,gg')       & \mapsto           &  (0,\ldots, g,\ldots ,0) \ , 
\end{array}
\end{equation}
where the unique nonzero entry $g$ is in the $g'$-th position. 
We also identify 
\begin{equation}
\label{eqn:6}
\begin{array}{ccc}
\Z[G\times G\times G] & \simeq & \bigoplus_{(g',g'')\in G\times G}\Z[G] \ , \\
(g,gg', gg'g'') & \mapsto  & (0,\ldots,g,\ldots, 0) \ , 
\end{array}
\end{equation}
where the unique nonzero   
entry $g$ is in the $(g',g'')$-th position. 

Let $M$ be a $G$-module.  Now the complex
$\Hom(\cC_{\bullet}^G,M)$ is identified with 
\begin{equation}
\label{eqn:2}
\cC^{\bullet}_{G,M}:= 
M\stackrel{d^0}{\lra}\bigoplus_{g'\in G}M\stackrel{d^1}{\lra} 
\bigoplus_{(g',g'')\in G\times G} M \cdots \ .
\end{equation}
Here the $g'$-th coordinate of the map $d^0$ is $m\mapsto g'\cdot m-m$
and the $(g',g'')$-th coordinate of $d^1$ is 
$(\ldots,m_g,\ldots)\mapsto g'\cdot m_{g''}-m_{g'g''}+m_{g'}$.
Of course, 
$H^i(G,M)$ is identified with the $i$-th cohomology of 
\eqref{eqn:2}. For instance, the kernel of $d^0$ is the module 
$M^{G}$ of $G$-invariants of $M$. 

Now let $H$ be a subgroup of $G$. Since restriction is an exact functor,
$\cC_{\bullet}^G$ is a resolution of $\Z$ as an $H$-module. 
We choose a set $Q\subset G$ 
of coset representatives, so $G=\bigcup_{q\in Q} Hq$. 

We have an isomorphism of $H$-modules
\begin{equation}
\label{eqn:hh}
\begin{array}{ccc}
\Z[G] & \simeq  &  \bigoplus_{q\in Q} \Z[H] \ ,  \\
hq    & \mapsto &  (0,\ldots , h,\ldots , 0) \ ,
\end{array}
\end{equation} 
where $h$ appears in the $q$-th
position ($h\in H,q\in Q$).  
Also 
\begin{equation}
\label{eqn:gg}
\begin{array}{ccc}
\Z[G\times G] & \simeq & \bigoplus_{(q,h',q')\in Q\times H\times Q} \Z[H] \ , \\
(hq, hh'q') & \mapsto & (0,\ldots , h,\ldots , 0) \ ,
\end{array}
\end{equation}
where $h$ appears in the $(q,h',q')$ position. 
We can project the resolution $\cC_{\bullet}^G$ 
to the standard resolution $\cC_{\bullet}^H$. 
Under the identification \eqref{eqn:hh} the map on the degree zero 
component is the sum of the $|Q|$ projection maps, and
under the identifications \eqref{eqn:4} and \eqref{eqn:gg}
the map on the degree 1 component sends 
the element 
$(0,\ldots , h,\ldots , 0)$ from \eqref{eqn:gg} to 
$(0,\ldots , h,\ldots , 0)$ with $h$ in the $h'$ position. 
Applying $\Hom_H(- , M)$ we get an inclusion of complexes
$\cC^{\bullet}_{H,M}$ into $\Hom_H(\cC_{\bullet}^G,M)$, and 
via our identifications, 
\begin{equation}
\begin{split}
\label{eqn:dd}
\xymatrix{ 
M\ar[r]\ar[d]^{\chi^0} & \bigoplus_{H} M \ar[r]\ar[d]^{\chi^1} &   \ldots \\
{\bigoplus_{Q} M} \ar[r] &  
{\bigoplus_{Q\times H\times Q} M} \ar[r]  & \ldots  \ .
}
\end{split}
\end{equation}
This allows us to take elements of $H^i(H,M)$, represented
as cocycles via the standard resolution, and realise them 
as cocycles in the complex $\Hom_H(\cC_{\bullet}^G,M)$.

Now we discuss cohomology of group extensions.
Assume that there is an exact sequence of groups
\begin{equation}
\label{eqn:HGQ}
1\ra H\ra G\ra Q\ra 1 \ .
\end{equation}
Then $Q$ acts on 
the cohomology $H^q(H,M)$ for all $q$, and 
there is an associated standard spectral sequence 
\begin{equation}
\label{eqn:sps}
E_2^{p,q}= H^p(Q,H^q(H,M))\Rightarrow H^{p+q}(G,M) \ .
\end{equation}
This leads to a $5$-term exact sequence 
\begin{equation}
\label{eqn:5}
0\ra H^1(Q,M^{H})\ra H^1(G,M)\ra H^1(H,M)^{Q}\stackrel{d_2^{0,1}}{\lra} 
H^2(Q,M^{H})\ra H^2(G,M) \ .  
\end{equation} 

The following is standard, but we will later use the formulas that are
explicitly given.

\begin{prop}
\label{prop:HGQ}
Given the exact sequence of finite groups \eqref{eqn:HGQ} and a
$G$-module $M$,
then
$$\tilde q\colon (\varphi\colon \Z[G^n]\to M)\mapsto
(\tilde q\cdot)\circ\varphi\circ(\tilde q^{-1}\cdot)$$
defines an action of $Q$ on the complex
$\Hom_H(\cC_{\bullet}^G,M)$, with invariants
$\Hom_G(\cC_{\bullet}^G,M)$.
When $G$ is a semi-direct product of $H$ and $Q$,
and we use \eqref{eqn:hh} and \eqref{eqn:gg} to identify
$\Hom_H(\cC_{\bullet}^G,M)$ with the
bottom row of \eqref{eqn:dd}, the action of $\tilde q\in Q$ is
given explicitly by
\begin{gather*}
\tilde{q}\cdot (\ldots, m_q, \ldots )  = (\ldots, \tilde{q}\cdot 
m_{\tilde{q}^{-1} q}, \ldots )\ ,  \\
\tilde{q}\cdot ( \ldots, m_{q,h',q'},\ldots ) = ( \ldots,
\tilde{q}\cdot m_{\tilde{q}^{-1}q, \tilde{q}^{-1}h'\tilde{q}, 
\tilde{q}^{-1}q'}, \ldots )\ .
\end{gather*}
\end{prop}

For many groups $G$ there are more efficient resolutions than
the standard resolution.
These are well known for finite abelian groups
(for instance, the case of bicyclic groups enters the calculations of
\cite{cubdiag}).
The following proposition captures the essential data needed to
compute $H^1$ and express elements there as cocycles for the
standard resolution, in the case of abelian groups with up to
three generators as well as dihedral groups.

\begin{nota}
\label{nota:g}
Let $G$ be a finite abelian 
group and $g\in G$ an element of order $n$.
Put $N_g:=1+g+\cdots + g^{n-1}$ and $\Delta_g := 1-g$ in $\Z[G]$.
For $g_1,\ldots,g_{\nu}\in G$ and $i_1, \ldots , i_{\nu}\in \Z$ 
the element in $\cC^G_1$ which, under the identification \eqref{eqn:4} 
is the vector $(0,\ldots, 1, \ldots, 0)$ with $1$ in the
$(g_{1}^{i_1}g_2^{i_2} \cdots g_{\nu}^{i_{\nu}})$-th  position, 
is denoted $\al_{i_1,\ldots, i_{\nu}}$. Similarly,
given $i'_1,\ldots , i'_{\nu}\in \Z$ 
the element in $\cC^G_2$ which, under the identification \eqref{eqn:6}
is the vector $(0,\ldots, 1, \ldots, 0)$ with $1$ in the
$(g_{1}^{i_1}g_2^{i_2} \cdots g_{\nu}^{i_{\nu}}, 
g_{1}^{i'_1}g_2^{i'_2} \cdots g_{\nu}^{i'_{\nu}})$-th position
is denoted $\al_{i_1,\ldots, i_{\nu},i'_1,\ldots, i'_{\nu}}$.
\end{nota}

\begin{prop}
\label{prop:efficient}
For each of the following classes of groups $G$ there exists a
resolution of $\Z$ by free $\Z[G]$-modules as stated.
In each case there is a morphism of complexes as indicated to this resolution
from the standard resolution.
\begin{itemize}
\item[(i)]
$G=\Z/n$, generated by $g\in G$:
$$
\cC^{[n]}_{\bullet} : =\cdots \Z[G]\stackrel{N_g}{\lra} \Z[G]
\stackrel{\Delta_g}{\lra}\Z[G] \ , 
$$
with
$\sigma_{\bullet}^{[n]}\,:\, \cC_{\bullet}^{G}\ra \cC_{\bullet}^{[n]}$
given by 
\begin{align*}
\sigma_1^{[n]}(\al_i) & = -1-g-\cdots - g^{i-1} \ , \\
\sigma_2^{[n]}(\al_{i,i'}) & = \left\{ 
\begin{array}{cc}  -1 & \text{ if } i+i'\ge n\ , \\
                   0 & \text {otherwise }\ . \end{array} \right. 
\end{align*}
\item[(ii)]
$G=\Z/n\oplus \Z/m$, with factors generated by $g$ and $h$:
$$
\cC^{[n,m]}_{\bullet} : =\cdots \Z[G]^3\stackrel{A^{[g,h]}}
\vvvlra \Z[G]^2
\stackrel{(\Delta_g\,\,\Delta_h)}\vvvlra\Z[G] \ , 
$$
where
$$
A^{[g,h]}:=\begin{pmatrix}
N_g&\Delta_h&0\\
0&-\Delta_g&N_h\end{pmatrix}
$$
with
$\sigma_{\bullet}^{[n,m]}\,:\, \cC_{\bullet}^{G}\ra \cC_{\bullet}^{[n,m]}$
given by 
$$
\sigma_1^{[n,m]}(\al_{i,j}) =
(-1-g-\cdots - g^{i-1}, -g^i(1+h+\cdots +h^{j-1})) \ .
$$
\item[(iii)]
$G=\Z/n\oplus \Z/m\oplus \Z/\ell$, with factors generated by
$g$, $h$, and $u$:
$$
\cC^{[n,m,\ell]}_{\bullet} : =\cdots \Z[G]^6\stackrel{
A^{[g,h,u]}}
\vvvvlra \Z[G]^3 
\stackrel{(\Delta_g  \,\, \Delta_h \,\,  \Delta_u)}{\vvvvvlra}
\Z[G] \ , 
$$
where
$$
A^{[g,h,u]}:=\begin{pmatrix}
N_g&\Delta_h&0&\Delta_u&0&0\\
0&-\Delta_g&N_h&0&\Delta_u&0\\
0&0&0&-\Delta_g&-\Delta_h&N_u\end{pmatrix}
$$
with
$\sigma_{\bullet}^{[n,m,\ell]}\,:\,
\cC_{\bullet}^{G}\ra \cC_{\bullet}^{[n,m,\ell]}$
given by 
$$
\sigma_1^{[n,m,\ell]}(\al_{i,j,k}) =
(-1-\cdots - g^{i-1}, -g^i(1+\cdots +h^{j-1}),
-g^ih^j(1+\cdots+u^k)) \ .
$$
\item[(iv)] $G={\mathfrak D}_n$, the dihedral group
generated by $g$ and $h$, with $g^n=h^2=(gh)^2=e$:
$$\cC^{\mathit{dih}[n]}_{\bullet} := \cdots \Z[G]^4\stackrel{D_n^3}\lra \Z[G]^3
\stackrel{D_n^2}\lra \Z[G]^2\stackrel{D_n^1}\lra \Z[G]\ ,$$
with
$$
D_n^3=\begin{pmatrix}
\Delta_g &  0 & 0 & N_h \\
0 & \Delta_h & 0 & -N_g \\
0 & 0 & \Delta_{gh} & -N_g\end{pmatrix},\qquad
D_n^2=\begin{pmatrix}
N_g & 0 & N_{gh} \\
0 & N_h & -N_{gh} \end{pmatrix},
$$
and $D_n^1=(\Delta_g\,\,\Delta_h)$, and
$\sigma_{\bullet}^{\mathit{dih}[n]}\,:\,
\cC_{\bullet}^{G}\ra \cC_{\bullet}^{\mathit{dih}[n]}$
given by
\begin{align*}
\sigma_1^{\mathit{dih}[n]}(\alpha_i) &=
(-1-g-\cdots - g^{i-1},0)\ ,\\
\sigma_1^{\mathit{dih}[n]}(\beta_i) &=
(-1-g-\cdots - g^{i-1},-g^i)\ ,
\end{align*}
where $\alpha_i$ is as in Notation \ref{nota:g} for the cyclic subgroup
generated by $g$, and where $\beta_i$ is the element of $\cC^G_1$
corresponding to $g^ih\in G$.
\end{itemize}
\end{prop}

\begin{proof}
All that is involved is checking, in each case, that we have
indeed specified (the tail end of) a resolution of $\Z$ as a
$\Z[G]$-module, and that the morphism from $\cC^G_{\bullet}$ is a
morphism of complexes.
\end{proof}

In each case ($G$ abelian or dihedral),
if we are given a $G$-module $M$, then applying $\Hom_G(-,M)$ to
the complex presented above
gives a practical method for computing group cohomology of $G$.
For instance, if $M$ is a $G$-module with $G=\Z/n$, generated by $g$,
then $H^i(G,M)$ is the $i{\rm th}$ cohomology
of
\begin{equation}\label{eqn:cohocyclic}
0\lra M\stackrel{\Delta_g}\lra M\stackrel{N_g}\lra M\lra \cdots.
\end{equation}

\begin{nota}
\label{nota:dual}
In the complex obtained by applying $\Hom_G(-,M)$, the maps
will be denoted as in Proposition \ref{prop:efficient}, but with the
super- and subscripts interchanged.
For example, $A_{[g,h]}\,:\,  M^2\ra M^3$ will denote the map
that sends the element $(m,0)$ to  
$(m+g\cdot m+\cdots+g^{n-1}\cdot m, m-h\cdot m, 0)$.
\end{nota}

By applying the efficient resolutions of Proposition \ref{prop:efficient}
to the group $Q$ acting on $\Hom_H(\cC_\bullet^G,M)$ in
Proposition \ref{prop:HGQ},
we can write the spectral sequence \eqref{eqn:sps} at the $E_0$ level.
This is necessary for computing $d_2^{0,1}$ in \eqref{eqn:5}, hence
for computing $H^1(G,M)$.

For example, when $Q$ is bicyclic we have:
\begin{coro}
\label{coro:E0}
If we have an extension of finite groups \eqref{eqn:HGQ} with $Q$ bicyclic,
then \eqref{eqn:sps} is the
spectral sequence of the bicomplex
$$
\xymatrix{ 
{\Hom_H(\Z[G^3],M)}       \ar[r] & 
{\Hom_H(\Z[G^3],M)^2} \ar[r] &  \cdots  \\
{\Hom_H(\Z[G^2],M)}       \ar[u]^{d_0^{0,1}}  \ar[r] &  
{\Hom_H(\Z[G^2],M)^2}\ar[u]^{d_0^{1,1}} \ar[r] & \cdots\ar[u]    \\
{\Hom_H(\Z[G],M)}               \ar[u]^{d_0^{0,0}}        \ar[r] &  
{\Hom_H(\Z[G],M)^2} \ar[u]^{d_0^{1,0}}  \ar[r] & 
{\Hom_H(\Z[G],M)^3}   \ar[u]
}     
$$
\end{coro}

\section{Computation of $\Br(S)/\Br(\Q)$ in the generic case}
\label{sect:comp-generic}
In this section we explain the computation of $H^1(G,M)$, where 
$M=\Pic(S_F)$, 
in the generic case $G=G_0$. 
We start by constructing, for each generator of $G$, the $8\times 8$ matrix
representing its action on $M$, referring to
Propositions \ref{prop:PicS} and \ref{prop:ga}.
In principle,  $H^1(G,M)$ can be computed using 
the standard resolution \eqref{eqn:2}.  In this case the map $d_1$
would be given by a $131072\times 1024$-matrix, which makes direct computations
impractical.  However, $G$ fits into a split exact sequence
\begin{equation}
\label{eqn:HGQ2}
1\to H\to G\to Q\to 1
\end{equation}
with $H=(\Z/4)^2\oplus(\Z/2)$ generated by
$\iota_a$, $\iota_b$, and $\iota_a\iota_b\iota_c$, and $Q=(\Z/2)^2$,
generated by $\sigma$ and $\tau$.
The technique of Section~\ref{sect:gr} simplifies the computation
considerably. 

\begin{prop}
\label{prop:genericH1}
For the generic Galois group $G=G_0$,
the cohomology group $H^1(G,M)$ is isomorphic to $\Z/2$.
\end{prop}

\begin{proof}
We use the $5$-term exact sequence \eqref{eqn:5}.
First we compute $M^H=M^G=\Z$, spanned by the anticanonical class. 
In particular, $H^1(Q,M^H)=0$.  Thus 
$H^1(G,M)$ is equal to the kernel of the map 
$$
d_2^{0,1}\,:\, H^1(H,M)^Q\ra H^2(Q,M^H) \ .
$$

We consider the diagram in Figure~\ref{fig:big},
where the bicomplex $E_0^{p,q}$ of Corollary \ref{coro:E0} is written using
the identifications \eqref{eqn:hh} and \eqref{eqn:gg}.
The group $H^1(H,M)$ is computed by the complex on the left side
of the diagram. In this diagram the horizontal arrows labeled 
$\sigma^i_{[4,4,2]}$ and $\chi^i$ give quasi-isomorphisms of complexes. 
The linear algebra required to compute 
$\Ker(M^3\ra M^6)$ is quite modest and the cohomology group is identified
as 
$$
H^1(H,M)= \Z/2 \ .
$$

It remains to take a single cocycle representative of the nonzero element of
$H^1(H,M)$
(necessarily $Q$-invariant in this case, though as noted below,
$Q$-invariance is tested a bit further on in the diagram chase)
and follow it through the diagram to determine whether it lies in 
the kernel of $d_2^{0,1}$. 

\begin{figure}
$$
\xymatrix{
M^6\ar[r]^(0.4){\sigma^2_{[4,4,2]}} & {\bigoplus_{H\times H} M}    \ar[r]  
& E_0^{0,2} \ar[r] &  E_0^{1,2} &   \\
M^3\ar[u]^{A_{[4,4,2]}}\ar[r]^(0.4){\sigma^1_{[4,4,2]}} 
& {\bigoplus_{H} M} \ar[u]\ar[r]^(0.4){\chi^1}    
& {\bigoplus_{Q\times H\times Q}  M} \ar[u]^{d_0^{0,1}}  \ar[r] &  
(\bigoplus_{Q\times H\times Q} M)^2 \ar[u]^{d_0^{1,1}} 
\ar[r] & E_0^{2,1}\\
M \ar[u] \ar@{=}[r]& M\ar[u] \ar[r]^(0.45){\chi^0}  & {\bigoplus_{Q}M}  
\ar[u]^{d_0^{0,0}}        \ar[r] &   
 (\bigoplus_{Q} M)^2 \ar[u]^{d_0^{1,0}}  \ar[r] & 
 (\bigoplus_Q M)^3 \ar[u]^{d_0^{2,0}} \\
         &              &     M^H  \ar@{^{(}->}[u]^{i_0} \ar[r] &
(M^H)^2 \ar@{^{(}->}[u]^{i_1} \ar[r]^{A_{[2,2]}} &
(M^H)^3 \ar@{^{(}->}[u]^{i_2}
}     
$$
\caption{$E_0$ spectral sequence}
\label{fig:big}
\end{figure}

We start with a representative in $M^3$ for the nontrivial element
$\lambda\in H^1(H,M)$, for instance,
$$   u=( ( 0,  0,  0,  0, -1, -1, -1,  1),
    ( 0,  0,  0,  0, -1,  1,  0,  0),
    ( 0,  0,  0,  0, -2,  0, -1,  1) ).$$
Let $v$ denote the image in $E_0^{1,1}$ of $u$ by the composite of
three horizontal maps in Figure \ref{fig:big}.
Now $v$ will in general lie in the image of $d_0^{1,0}$ if and only if
$\lambda$ is $Q$-invariant.
In this case, a linear algebra solver produces
$$   v_0=( ( 0, 0, 0, 0, -1, 1, 0, 0 )^{*4} ,
           ( 0, 0, 0, 0, -1, -1, -1, 1 ) ^{*4} ) )$$
satisfying $d_0^{1,0}(v_0)=v$, where each vector with superscript $*4$
denotes the element in $\bigoplus_QM$ with the vector repeated 4 times.
Applying the cobounday map $E_0^{1,0}\to E_0^{2,0}$ to $v_0$
necessarily produces an element in the image of $i_2$,
representing $d_2^{0,1}(\lambda)$ in $H^2(Q,M^H)$.
This can be tested for being a coboundary; in the present case we get
$0$ exactly.
So $d_2^{0,1}$ is trivial, and $H^1(G,M)=\Z/2$.
\end{proof}

\begin{coro}
\label{coro:genericS}
If the del Pezzo surface $S$ given by \eqref{eqn:d2} is general,
meaning that the hypotheses of Proposition \ref{prop:ga} are met, then we have
$$\Br(S)/\Br(\Q)=\Z/2.$$
\end{coro}

In Example \ref{exam:obstr3}, below, we will see how to construct explicitly
an Azumaya algebra representing the nontrivial element of
$\Br(S)/\Br(\Q)$ and use it to test the Brauer-Manin obstruction.

\section{The non-generic case}
\label{sect:non-gen}
We start by presenting some
examples when the Galois group is smaller than in the generic case.

\begin{exam}
\label{exam:n-1}
Consider the case
$(A,B,C)=(-6,-3,2)$.
The Galois group of the field $F$, defined in \eqref{eqn:F},
has order 32; it is an extension of the Klein four-group by
$(\Z/4)\oplus (\Z/2)$. It is possible to write $G$ as a split extension
$$
 1 \ra H \ra G \ra \Z/2 \ra 1
$$
where $H=(\Z/4)^2$, generated by 
$\iota_a\iota_b$ and $\sigma\tau\iota_a\iota_c$,
and $\Z/2$ is 
generated by $\sigma$. 
In this case, we compute $H^1(H,M)=0$. By \eqref{eqn:5},  
$H^1(G,M)$ is isomorphic to $H^1(\Z/2, M^H)$.
We find that $M^H$ has rank 2, spanned by
$$
( -1,  -1,  -1,  -1,  -1,  -1,  -1, 3),\qquad ( 1,  1,  1,  1,  1,  1,  0, -2),
$$
hence $M^H$ is isomorphic to $\Z\oplus \Z'$, where $\Z'$ is
free of rank 1 with nontrivial $\Z/2$-action.
So, we have
$$H^1(G,M)=\Z/2.$$ 
As in the generic case, we have $M^G=\Z$, that is,
$\Pic(S)$ has rank 1.
\end{exam}

\begin{exam}
\label{exam:n-2}
The case $(A,B,C)=(1,1,-2)$ is interesting because
$\Pic(S)$ has rank 2.
The Galois group $G$ fits into an exact sequence
$$1\to\Z/4\to G\to \Z/2\to 1$$
with subgroup $H=\Z/4$ generated by $\iota_c\sigma\tau$
and $\Z/2$ generated by $\tau$.
As in Example~\ref{exam:n-1} we have
$H^1(H,M)=0$.
Now $M^H$ has rank 3, with generators
$$
( -1,  -1,  -1,  -1,  -1,  -1,  -1, 3)\ ,\quad
( 0,  0,  0,  0,  1, -1,  0,  0)\ ,\quad
( 0,  0,  0,  0,  1,  1,  1, -1)\ ,
$$
and the action of $\tau$ fixes the first 2 vectors and
negates the third.
Hence
$$H^1(G,M)=\Z/2,$$
and $\Pic(S)$ has rank 2.
\end{exam}

\begin{exam}
\label{exam:n-3}
The case $(A,B,C)=(1,1,1)$ yields $G=\Gal(\Q(\zeta)/\Q)$,
the Klein four-group, and we directly compute
$$H^1(G,M)=(\Z/2)^3.$$
In this case $\Pic(S)$ has rank 1.
\end{exam}

The comprehensive treatment proceeds via a case-by-case computer analysis
of subgroups of the generic Galois group.
We obtain the following, as our main result.
\begin{theo}
\label{theo:big}
Let $S$ have the form \eqref{eqn:d2}, where $A$, $B$, and $C$ are
nonzero integers.
Then $\Br(S)/\Br(\Q)$ is isomorphic to one of the following groups:
\begin{gather*}
(1), \ \Z/2, \ \Z/4, \ (\Z/2)\oplus(\Z/2), \\
(\Z/4)\oplus(\Z/2), \ (\Z/2)\oplus(\Z/2)\oplus(\Z/2).
\end{gather*}
Also, $\Br(S)/\Br(\Q)$ is nontrivial in every case where
$\Pic(S)$ is isomorphic to $\Z$.
\end{theo}

\begin{proof}
Given such $S$, a choice of $4$-th roots $a$, $b$, $c$ of the
coefficients leads to a realisation of
$G=\Gal(\Q(\zeta,a^2,b/a,c/a)/\Q)$ as a subgroup of the
generic Galois group $G_0$.
This will be a subgroup mapping surjectively to $Q$ via the map
in \eqref{eqn:HGQ2}.
Computer analysis reveals that
every subgroup of $G_0$ which maps surjectively to $Q$ can be expressed as a
semi-direct product of abelian groups.

We recognise that the same group cohomology
must arise from any two subgroups which
differ by conjugation, or by the obvious outer action of ${\mathfrak S}_3$ on $G_0$
corresponding to permutations of the
$x$, $y$, and $z$ coordinates.
More generally, the full group of automorphisms of $\Pic(S_{\overline\Q})$
preserving the intersection pairing and the anticanonical class is
the Weyl group $W(E_7)$; see \cite{kst}.
This is a group of order $2903040$,
generated by $G_0$ and the group ${\mathfrak S}_7$ of permutations
of $v_1$ through $v_7$.
Any two subgroups of $G_0$ that are conjugate in $W(E_7)$ must have
same cohomology.
There are $194$ classes of subgroups of $G_0$, up to conjugation in
$W(E_7)$, which contain a group that surjects onto $Q$.
When the methods of section \ref{sect:gr} are applied to a representative
of each class of subgroups, the $H^1$ group that results
is always one of the groups listed in the statement of the theorem.
Moreover, the trivial group arises as $H^1(G,M)$ only in cases with the
rank of $M^G$ greater than or equal to $2$.

There are too many classes of subgroups to list them all,
so we content ourselves with displaying,
in Table \ref{tab:maxsubG}, all the maximal subgroups of $G_0$
that surject onto $Q$, up to the ${\mathfrak S}_3$-action.
These are grouped by conjugacy in $W(E_7)$.
For each subgroup,
we display the cohomological invariants,
the condition on $A$, $B$, and $C$ that forces the Galois group to be
contained in the subgroup,
and a representative $(A,B,C)$ for this subgroup.
The complete list of subgroups,
together with accompanying \texttt{magma} code,
can be found under the computing link
at the first author's web page
\texttt{http://www.maths.warwick.ac.uk/\~{}kresch/}.
\end{proof}

\begin{table}
$$
\begin{array}{c|c|c|c|c}
G & \Br(S)/\Br(\Q) & \Pic(S) & \mathrm{Condition} & \mathrm{Example} \\
\hline
\langle
\iota_a\sigma, \iota_a\iota_b, \iota_a\iota_c, \tau
\rangle &
\Z/2 &
 \Z &
-2ABC\in(\Q^*)^2 &
(-15,10,3) \\
\langle
\iota_a\sigma, \iota_b, \iota_c, \tau
\rangle &
\Z/2 &
 \Z &
-2A\in(\Q^*)^2 &
(-2,3,5) \\
\hline
\langle
\iota_a\sigma, \iota_a\iota_b, \iota_c, \tau
\rangle &
\Z/2 &
 \Z &
-2AB\in(\Q^*)^2 &
(-6,3,5) \\
\hline
\langle
\iota_a\tau, \iota_a\iota_b, \iota_a\iota_c, \sigma
\rangle &
\Z/2 &
 \Z &
2ABC\in(\Q^*)^2 &
(3,10,15) \\
\langle
\iota_a\tau,\iota_b,\iota_c,\sigma
\rangle &
\Z/2 &
 \Z &
2A\in(\Q^*)^2 &
(2,3,5) \\
\hline
\langle
\iota_a\tau, \iota_a\iota_b, \iota_c, \sigma
\rangle &
\Z/2 &
 \Z &
2AB\in(\Q^*)^2 &
(-6,-3,5) \\
\hline
\langle
\iota_a\sigma, \iota_a\iota_b, \iota_a\iota_c, \sigma\tau
\rangle &
(1) &
 \Z\oplus\Z &
-ABC\in(\Q^*)^2 &
(-15,3,5) \\
\hline
\langle
\iota_a\sigma, \iota_b, \iota_c, \sigma\tau
\rangle &
\Z/2 &
 \Z &
-A\in(\Q^*)^2 &
(-1,3,5) \\
\hline
\langle
\iota_a\sigma, \iota_a^2, \iota_a\iota_b, \iota_c, \sigma\tau
\rangle &
\Z/2 &
 \Z &
-AB\in(\Q^*)^2 &
(-63,7,15) \\
\hline
\langle
\iota_a\iota_b, \iota_a\iota_c, \sigma, \tau
\rangle &
\Z/2 &
 \Z &
ABC\in(\Q^*)^2 &
(3,5,15) \\
\hline
\langle
\iota_b, \iota_c, \sigma, \tau
\rangle &
\Z/2 &
 \Z &
A\in(\Q^*)^2 &
(1,3,5) \\
\hline
\langle
\iota_a^2, \iota_a\iota_b, \iota_c, \sigma, \tau
\rangle &
\Z/2\oplus\Z/2 &
 \Z &
AB\in(\Q^*)^2 &
(-63,-7,5) \\
\hline
\end{array}
$$
\caption{Possible Galois groups among maximal subgroups of $G_0$}
\label{tab:maxsubG}
\end{table}

\begin{rema}
The significance of $\Pic(S)$ being isomorphic to $\Z$, according to
the Enriques-Manin-Iskovskih
classification of surfaces, is that these are the minimal surfaces which
are not conic bundles.
\end{rema}

\section{Examples of Brauer-Manin obstruction}
\label{sec:examples}
Here we compute the Brauer-Manin obstruction to the
Hasse principle in several representative cases.

\begin{exam}
\label{exam:obstr1}
The case $(A,B,C)=(-25,-5,45)$.
The group $G=\Gal(F/\Q)$ has order 32 and fits into an
exact sequence
$$
1\to H\to G\to \Z/2\to 1
$$
with $H=(\Z/4)\oplus(\Z/2)^2$, generated by $\iota_a^2\iota_b\iota_c$,
$\iota_c^2$, and $\sigma\tau$, and $\Z/2$ generated by $\sigma\iota_a\iota_b$.
Computing, as in the previous section, we find
\begin{equation}
\label{firstmapiso}
H^1(\Z/2, M^H)\eqto H^1(G,M)
\end{equation}
in the sequence \eqref{eqn:5}, with $M^H$ equal to the span of
$(-1,-1,-1,-1,-1,-1,-1,3)$ and $(1,1,1,1,1,1,0,-2)$.
Hence, as in Example \ref{exam:n-1}, we have
$$H^1(G,M)=\Z/2.$$

Because of \eqref{firstmapiso}, there will exist
a class in $\Br(S)$, not in $\Br(\Q)$, which is
annihilated by the field extension
$\Q\to\Q[i]=F^H$.
This makes it convenient to carry out the procedure described
in \cite[Chap.\ VI]{manin} for constructing
a central simple algebra over the function field of $S$ which
is the restriction of a sheaf of Azumaya algebras that is
nontrivial in $\Br(S)/\Br(\Q)$.
This works as follows.
Let $\alpha$ be any divisor on $S_{\Q[i]}$
whose class in $M=\Pic(S_F)$ is $(1,1,1,1,1,1,0,-2)$.
Since $\alpha$ and its complex conjugate
$\overline{\alpha}$ sum to $0$ in $M$ there will exist
a rational function $g$ whose divisor is $\alpha+\overline{\alpha}$;
then we consider the quaternion algebra
$$(-1,g)\in\Br(S).$$

We can take $\alpha$ to be the class of a conic minus an anticanonical
divisor,
$$\alpha=D-(z=0),$$
where the conic $D\subset S_{\Q[i]}$ is taken to lie above a conic meeting
the Fermat quartic in $4$ tangencies, e.g.,
$$-5x^2-2y^2+9z^2=0, \qquad\qquad w=i(3y^2-6z^2).$$
It is easy to check that $\alpha$ has the correct class in $M$.
With this choice, we can take
$$g=-5(x/z)^2-2(y/z)^2+9.$$

By the geometry underlying the choice of $D$, we have $g>0$ for any
$[x:y:z]\in\PP^2(\Q)$ that has real points over it in $S$.
It is now only necessary to complete $p$-adic analyses at the primes $p=2$ and
$p=3$
(since $5$-adically, $\Q[i]$ is a split extension of $\Q$).
For the $2$-adic analysis, we assume $x$, $y$, and $z$ to be $2$-adic
integers, not all even,
and find by analysis mod $16$ that the condition
$-25x^4-5y^4+45z^4$ should be a $2$-adic square implies
$x$ and $z$ are odd and $y$ is even.
So, without loss of generality, we may take $z=1$.
By mod $32$ analysis, the only possible values of $(x,y)$ mod $8$ are
$$
\begin{array}{cccccccc}
(1,2),&(1,6),&(3,0),&(3,4),&
(5,0),&(5,4),&(7,2),&(7,6).
\end{array}
$$
In each case we find $g=12 \pmod{16}$, hence $(-1,g)$ is
{\em ramified} at all $2$-adic points of $S$.
By a similar analysis mod 27 we find that at any $3$-adic point
$x$ and $y$ are prime to $3$, hence so is $g$, and
$(-1,g)$ is {\em unramified} at all $3$-adic points of $S$.
Therefore $S$ provides an example of Brauer-Manin obstruction to the
Hasse principle.
\end{exam}

\begin{exam}
\label{exam:obstr2}
Here we show that Example \ref{exam:n-1} fits into an infinite
family of examples of Brauer-Manin obstruction to the Hasse principle.
Consider
$$(A,B,C)=(-2p,-p,2),$$
where $p$ is any prime such that
$$p=3\pmod{16}.$$
The computation of the group cohomology is exactly as in
Example \ref{exam:n-1}.
So, $H^1(G,M)=H^1(\Z/2,M^H)=\Z/2$.
We proceed as in Example \ref{exam:obstr1}.

By the condition on $p$ we may write
$$p=u^2+2v^2$$
for positive integers $u$ and $v$, necessarily both odd.
Define $s=(-1)^{(u-v)/2}$.
Solving for the plane conic tangent to the quartic at the
points $(\pm \sqrt{su/p},\pm\sqrt{2v/p})$, we find that with the curve $D$
given by
$$-sux^2-vy^2+z^2=0, \qquad\qquad w=i(-2vx^2+suy^2),$$
the cycle $D-(z=0)$ has class $(1,1,1,1,1,1,0,-2)$ in $M$.
Set
$$g=-su(x/z)^2-v(y/z)^2+1.$$
Then $(-1,g)$ is
\begin{itemize}
\item[(i)] {\em unramified} at real points of $S$;
\item[(ii)] {\em ramified} at all $2$-adic points of $S$;
\item[(iii)] {\em unramified} at all $p$-adic points of $S$;
\end{itemize}
and there is a Brauer-Manin obstruction to the Hasse principle.

We leave the verification of (i)--(ii) to the reader.
For (iii) we need the following lemma.
\begin{lemm}\label{ntlemma}
Let $p$ be a prime with $p=3\pmod{16}$.
Write
$p=u^2+2v^2$ for positive integers $u$ and $v$.
Now, if we let $y$ be a solution to $y^4=-2\pmod{p}$ then we have
$vy^2=(-1)^{(u-v)/2}u\pmod{p}$.
\end{lemm}

\begin{proof}
The two square roots of $-2$ mod $p$ are $\pm uv^{-1}$.
So $y^2=\pm uv^{-1}\pmod{p}$ and the lemma is asserting that
the correct sign is $(-1)^{(u-v)/2}$,
or equivalently, that
\begin{equation}
\label{eqn:uv}
\Bigl(\frac{uv}{p}\Bigr)=(-1)^{(u-v)/2}.
\end{equation}
By quadratic reciprocity,
$$\Bigl(\frac{u}{p}\Bigr)=(-1)^{(u-1)/2}\Bigl(\frac{p}{u}\Bigr)
\qquad\text{and}
\qquad\Bigl(\frac{v}{p}\Bigr)=(-1)^{(v-1)/2}\Bigl(\frac{p}{v}\Bigr).$$
If $p'$ is a prime dividing $v$, then $p$ is a quadratic residue mod $p'$.
This and a similar consideration when $p'$ divides $u$ yield
$$\Bigl(\frac{p}{v}\Bigr)=1\qquad\text{and}\qquad
\Bigl(\frac{2p}{u}\Bigr)=1.$$
By mod $16$ analysis, $u=\pm 1\pmod{8}$, hence
$\bigl(\frac{2}{u}\bigr)=1$.
So, \eqref{eqn:uv} holds.
\end{proof}

To establish (iii) we claim that
for any $p$-adic integer solution $(w,x,y,z)$ to
\eqref{eqn:d2},
with not all of $w$, $x$, $y$, and $z$ divisible by $p$, 
the $p$-adic integer $z^2g=-sux^2-vy^2+z^2$ is not divisible by $p$.
Indeed, since $2$ is not a quadratic residue mod $p$ we must have
$p$ dividing $z$,
hence $x$ and $y$ are nonzero mod $p$.
Without loss of generality we
suppose $x=1$. Now $y$ must be a $4$-th root of $-2$ mod $p$.
The claim follows from Lemma \ref{ntlemma}.
\end{exam}

\begin{exam}
\label{exam:obstr3}
Here we give a recipe for testing the presence of Brauer-Manin obstruction
to the Hasse principle
in the generic case, i.e., when the Galois group has order 128.
This occurs precisely when the set
$$\bigl\{\,A^\alpha B^\beta C^\gamma (-1)^\delta 2^\varepsilon\,\bigm|\,
(\alpha,\beta,\gamma,\delta,\varepsilon)\in \{0,1\}^5\smallsetminus
\{(0,0,0,0,0)\}\,\bigr\}$$
contains no perfect squares (see Table \ref{tab:maxsubG}).

Let $S$ be such a surface, and assume $S$ has rational points in all
completions of $\Q$.
By Corollary \ref{coro:genericS},
we have $\Br(S)/\Br(\Q)=H^1(G,M)=\Z/2$.
We use the fact that $G$ has a subgroup of index two
$$H=\langle \sigma\tau, \iota_a^2, \iota_a\iota_b, \iota_a\iota_c,
\iota_a\sigma\rangle$$
with the property that
$$
M^H=\langle\,(-1,-1,-1,-1,-1,-1,-1,3),(1,1,1,1,1,1,0,-2)\,\rangle,
$$
and hence
$H^1(G/H,M^H)\eqto H^1(G,M)$. Therefore, we can construct a
quaternion algebra as in Example \ref{exam:obstr1}.
In this case,
$$F^H=\Q(\sqrt{-ABC}).$$

Let $\theta=\sqrt{-ABC}$, and let $(r_0:s_0:t_0)$
be a $\Q(\theta)$-rational point on the conic
\begin{equation}
\label{eqn:rst}
Ar^2+Bs^2+Ct^2=0.
\end{equation}
By our assumption on $S$,
such a point exists by the Hasse principle:
local solutions to \eqref{eqn:rst} arise by rewriting \eqref{eqn:d2} as
$((\theta z^2)^2+ABw^2)/((By^2)^2+ABx^4)=A$.
Now
$$Ar_0x^2+Bs_0y^2+Ct_0z^2=0$$
defines a conic over $\Q(\theta)$, meeting the quartic curve
in tangencies.
By the identity
\begin{multline*}
C^2t_0^2(Ax^4+By^4+Cz^4)+ABC(s_0x^2-r_0y^2)^2 \\
+C(Ar_0x^2+Bs_0y^2+Ct_0z^2)(Ar_0x^2+Bs_0y^2-Ct_0z^2)=0,
\end{multline*}
there is a curve $D$ on $S_{\Q(\theta)}$ defined by
$$Ar_0x^2+Bs_0y^2+Ct_0z^2=0, \qquad\qquad
w=\theta (s_0x^2-r_0y^2)/(Ct_0)$$
such that the union of $D$ and its conjugate is rationally equivalent to
twice the anticanonical class.
This rational equivalence is given explicitly by the rational function
$$
g:=(Ar_1s_1+A^2BCr_2s_2)+(Bs_1^2-A^2BCr_2^2)(y/x)^2+
Cs_1t_0(z/x)^2+ACr_2t_0w/x^2,
$$
where we suppose $t_0\in\Q$ and write
$$r_0 = r_1+r_2\theta\qquad{\rm and}\qquad
s_0 = s_1+s_2\theta.$$
To test the Brauer-Manin obstruction to the Hasse principle
for $S$, one has to analyse the quaternion algebra
$$(-ABC,g)$$
at real- and $\Q_p$-valued points of $S$ (for $p$ dividing $2ABC$).

We give an example of
nontrivial Brauer-Manin obstruction in this case.
Consider
$(A,B,C)=(-126,-91,78)$.
Then we may take $r_0=-13$, $s_0=-12$, and $t_0=21$,
and $g$ is proportional to
$$3+2(y/x)^2+3(z/x)^2.$$
In this case the quaternion algebra $(-ABC,3+2(y/x)^2+3(z/x)^2)$
is ramified at all $\Q_2$-points of $S$ and
unramified at all points in all other completions.
\end{exam}

\begin{exam}
\label{exam:obstr4}
The case $(A,B,C)=(34,34,34)$.
Here $G=\Gal(F/\Q)$ is isomorphic to $(\Z/2)^3$:
$$G=\langle \iota_a \iota_b \iota_c \sigma, \tau, \sigma\rangle.$$
We have $H^1(G,M)=(\Z/2)^3$.
In fact, for the index-two subgroup
$$H=\langle \iota_a \iota_b \iota_c \sigma, \tau\rangle$$
we have
$M^H$ spanned by
\begin{equation}\label{rankfourbasis}
\begin{array}{c}
(1,-1,0,0,0,0,0,0),\\
(0,0,1,-1,0,0,0,0),\\
(0,0,0,0,1,-1,0,0),\\
(-1,-1,-1,-1,-1,-1,-1,3),
\end{array}
\end{equation}
and
$$H^1(G/H,M^H)\eqto H^1(G,M).$$
Here, $\sigma$ in $G/H$
acts nontrivially on the first three vectors
in \eqref{rankfourbasis} and trivially on the last.
We have
$$F^H=\Q(\sqrt{-17}).$$

Using \eqref{eqn:cohocyclic} we can identify elements of
$\Br(S)/\Br(\Q)$ with the image of the $(-1)$-eigenspace of $M^H$
(under the $\sigma$-action).
To produce quaternion algebras representing a given element of
$\Br(S)/\Br(\Q)$ we need to find divisors defined over $\Q(\sqrt{-17})$
representing particular classes in $M^H$.
Notice that the class of any combination of exceptional curves defined over
$\Q(\sqrt{-17})$ in $M^H$ is a coboundary of \eqref{eqn:cohocyclic}.
Hence, we need additional cycles defined over $\Q(\sqrt{-17})$.
We use descent to produce line bundles on $S_{\Q(\sqrt{-17})}$
and obtain the desired cycles
as loci of vanishing of rational sections of these line bundles.

Here we explicitly carry out the task of representing the
class of the first entry of \eqref{rankfourbasis} in $\Br(S)$.
Set $\rho=\iota_a \iota_b \iota_c \sigma$.
Over $F=\Q(\sqrt{-17},\zeta)$ we have
\begin{equation}
\label{eqn:itsclass}
[L_{x,\zeta,+}]-[L_{x,\zeta^3,-}]=(1,-1,0,0,0,0,0,0)
\end{equation}
in $\Pic(S_F)$.
Consider the
line bundle $\cO([L_{x,\zeta,+}]-[L_{x,\zeta^3,-}])$ together with
isomorphisms
$$\cO(L_{x,\zeta,+}-L_{x,\zeta^3,-})\stackrel{\eta}\lra
\cO(L_{x,\zeta^7,-}-L_{x,\zeta^5,+})$$
and
$$\cO(L_{x,\zeta,+}-L_{x,\zeta^3,-})\stackrel{\xi}\lra
\cO(L_{x,\zeta^3,+}-L_{x,\zeta,-}).$$
These constitute descent data
(for the covering $S_F\to S_{\Q(\sqrt{-17})}$)
provided that the diagram
$$\xymatrix{
{\cO(L_{x,\zeta,+}-L_{x,\zeta^3,-})}  \ar[r]_{\eta}  \ar[d]^{\xi}
\ar@/^18pt/@<2pt>[rr]^1 \ar@/_24pt/@<-35pt>[dd]_1 &
{\cO(L_{x,\zeta^7,-}-L_{x,\zeta^5,+})}\ar[r]_{\rho(\eta)}\ar[d]^{\rho(\xi)}&
{\cO(L_{x,\zeta,+}-L_{x,\zeta^3,-})} \\
{\cO(L_{x,\zeta^3,+}-L_{x,\zeta,-})} \ar[r]_{\tau(\eta)}\ar[d]^{\tau(\xi)} &
{\cO(L_{x,\zeta^5,-}-L_{x,\zeta^7,+})}
\\
{\cO(L_{x,\zeta,+}-L_{x,\zeta^3,-})}
}$$
commutes.
The isomorphisms given by
$$\eta=\delta
\frac{x^2-iy^2+z^2-\frac{1}{\sqrt{34}}w}{x^2-iy^2-z^2+\frac{1}{\sqrt{34}}w}
\qquad\text{and}\qquad
\xi=\varepsilon\frac{\zeta y + z}{\zeta^3 y + z},$$
satisfy this condition if and only if
$\delta$, $\eta\in F$ satisfy
\begin{align}
\delta \, \rho(\delta) &= -1, \label{eqn:ni1}\\
\varepsilon \, \tau(\varepsilon) &= 1, \label{eqn:ni2}\\
\delta\, \rho(\varepsilon) &= \tau(\delta)\,\varepsilon. \label{eqn:ni3}
\end{align}
One solution to \eqref{eqn:ni1}--\eqref{eqn:ni3} is
$$\delta = \sqrt{-17}\zeta - 4\zeta^3 \qquad\text{and}\qquad
\varepsilon = 4\zeta + \sqrt{-17}\zeta^3.$$
This yields, by effective descent, a line bundle
$\cE$ on $S_{\Q(\sqrt{-17})}$.

Using \eqref{eqn:ni1}--\eqref{eqn:ni3} and descent, we see that
$$f:=1+\rho(\eta)+\tau(\xi)+\rho(\eta\,\tau(\xi))$$
defines a rational section of $\cE$.
We write $f$ as a quotient of quartic polynomials and observe that
$f$ has (with respect to local trivializations of $\cE$)
a simple pole along
$L_{x,\zeta,-}\cup L_{x,\zeta^3,-}\cup L_{x,\zeta^5,+}\cup L_{x,\zeta^7,+}$
and a zero of order one along some curve $Z$.
Then, by \eqref{eqn:itsclass},
we deduce that
$$[Z]=(-3,-1,-2,-2,-2,-2,-2,6)$$
in the Picard group.
Therefore, if $h\in\Q(S)$ defines a rational equivalence between
$Z\cup \sigma(Z)$ and some hyperplane sections,
then the quaternion algebra $(-17,h)$ represents an element of
$\Br(S)$ of the desired class in $\Br(S)/\Br(\Q)$.

Denoting by $g$ the numerator of $f$, we have
\begin{multline*}
g=
(x^2+iy^2+z^2+\frac{1}{\sqrt{34}}w)
[y^2+iz^2+(4\zeta-\sqrt{-17}\zeta^3)(y^2+\sqrt{2}yz+z^2)]\\
+(x^2+iy^2-z^2-\frac{1}{\sqrt{34}}w)
[y^2+\sqrt{2}yz+z^2+(4\zeta-\sqrt{-17}\zeta^3)(-y^2+iz^2)].
\end{multline*}
The simultaneous vanishing of $g$, $\rho(g)$,
$\tau(g)$, and $\rho\tau(g)$ defines the curve $Z$.
Equivalently, writing
$$g=p_0 + p_1\zeta + p_2\zeta^2 + p_3\zeta^3$$
with $p_i\in \Q(\sqrt{-17})[w,x,y,z]$ we have $Z$ defined by the
vanishing of $p_i$ for $i=0$, $\ldots$, $3$.
A unique (up to scale) $\Q(\sqrt{-17})$ linear combination
of these is defined over $\Q$, namely
\begin{align*}
h_1&:=\frac{1}{2} p_0 + \frac{4-\sqrt{-17}}{2} p_1 +
\frac{1}{2} p_2 - \frac{4+\sqrt{-17}}{2} p_3 \\
&=wy^2+wz^2+x^2y^2+8x^2yz+x^2z^2+y^4-z^4.
\end{align*}
Then $h=h_1/x^4$ is as desired.
Cyclically permuting the variables $x$, $y$, and $z$, we obtain
polynomials $h_2$ and $h_3$ such that the classes of
$(-17,h_i/x^4)$ generate $\Br(S)/\Br(\Q)$.

The ramification pattern of an Azumaya algebra is an invariant of
its class in $\Br(S)$.
However, in practice, the ramification pattern of an algebra
$(-17,h_i/x^4)$ is difficult to test on $p$-adic points where
$h_i$ vanishes to high order.
Hence it is helpful to have multiple rational
functions determining the same class in $\Br(S)/\Br(\Q)$.
We can obtain additional functions by repeating the
previous construction for different
solutions to \eqref{eqn:ni1}--\eqref{eqn:ni3}.
For instance, $(-\delta,\varepsilon)$ is another solution.
If we carry out the above procedure with this solution we obtain
$$h_4 = wy^2+wz^2+x^2y^2+8x^2yz+x^2z^2-y^4+z^4,$$
with the property that
$(-17,h_1/x^4)$ and $(-17,h_4/x^4)$ are equal in $\Br(S)/\Br(\Q)$.
We obtain $h_5$ and $h_6$ similarly: the effect of the full set of
permutations of $x$, $y$, and $z$ is that we now have two representatives of
each of the generators of $\Br(S)/\Br(\Q)$.
We let ${\mathfrak q}_i\in \Br(S)$ denote
$(-17,h_i/x^4)$, for each $i$.

To gain full advantage of having these classes, we need to know how
${\mathfrak q}_i$ and ${\mathfrak q}_{i+3}$ differ in $\Br(S)$.
This is discovered by finding a relationship that makes explicit
their equality in $\Br(S)/\Br(\Q)$.
Using linear algebra, we have identified a rational equivalence on
$S_{\Q(\sqrt{-17})}$ between
$Z$ and the analogous curve for the function $h_4$; its norm relates
$h_1$ and $h_4$ modulo the defining equation of $S$:
\begin{align*}
h_1h_4&=\frac{1}{9}\Bigl[
\bigl(\frac{1}{2}wy^2+4wyz+\frac{1}{2}wz^2+17x^2y^2+17x^2z^2
-4y^4+y^3z+yz^3-4z^4\bigr)^2\\
&{}+17\bigl(\frac{1}{34}wy^2+\frac{4}{17}wyz+\frac{1}{34}wz^2+x^2y^2+x^2z^2
+4y^4-y^3z-yz^3+4z^4\bigr)^2\Bigr]\\
&{}+(-33y^4+16y^3z-2y^2z^2+16yz^3-33z^4)\bigl(x^4+y^4+z^4-\frac{1}{34}w^2\bigr).
\end{align*}
Similar identities hold under cyclic permutations of $x$, $y$, and $z$,
and we thus have
$${\mathfrak q}_i={\mathfrak q}_{i+3}$$
in $\Br(S)$, for each $i$.

Here are the results of the local analysis, confirming the
presence of a Brauer-Manin obstruction:
\begin{itemize}
\item ${\mathfrak q}_i$ is unramified on points of $S(\R)$ with $w>0$
and ramified on points with $w<0$, for all $i$.
\item $S(\Q_2)$ is the disjoint union of two nonempty sets,
$U$ and $R$, such that each ${\mathfrak q}_i$ is unramified
on $U$, and each ${\mathfrak q}_i$ is ramified on $R$.
\item At any point of $S(\Q_{17})$, exactly two of
$\{{\mathfrak q}_1,{\mathfrak q}_2,{\mathfrak q}_3\}$ are ramified.
\end{itemize}
\end{exam}

\begin{rema}
\label{rema:howtosolveni}
We produced the solution to \eqref{eqn:ni1}--\eqref{eqn:ni3} by inspection.
A more systematic way to proceed would be to solve just \eqref{eqn:ni1},
obtaining by descent a line bundle defined over $\Q(\sqrt{-17},\sqrt{2})$.
Descending further to $\Q(\sqrt{-17})$ then hinges upon solving a norm equation
for the quadratic extension $\Q(\sqrt{-17})\to \Q(\sqrt{-17},\sqrt{2})$.
\end{rema}

\begin{exam}\label{exam:obstr5}
The case $(A,B,C)=(-9826,-2,136)=(-2p^3,-2,8p)$ with $p=17$
illustrates working with a non-cyclic Azumaya algebra.
We have $F=\Q(\zeta,\sqrt[4]p)$.
The Galois group of $F$ over $\Q$ has order 16:
$$G=\langle \iota_a\iota_b\iota_c\sigma\tau, \iota_a^3\iota_c,
\iota_b\iota_c^3\sigma\rangle.$$
In this case, $H^1(G,M)=\Z/4$.
There is no Brauer-Manin obstruction coming from $2$-torsion in $\Br(S)$.
Indeed, the motivated reader can produce a subgroup $H$ of index $2$ in $G$
with $H^1(G/H,M^H)=\Z/2$ and show that
$(-2,136+(y/x)^2+18(z/x)^2)$ generates the $2$-torsion in $\Br(S)/\Br(\Q)$,
yet is unramified at all points $S$ in every completion of $\Q$.
This means that the obstruction analysis requires a representative of a
generator of $\Br(S)/\Br(\Q)$.

The central element $u:=\iota_a\iota_b\iota_c\sigma\tau$ of $G$ satisfies
$$F^u=\Q(i,\sqrt[4]p)$$
and
$$H^1(G/\langle u\rangle, M^u)=\Z/4.$$
We remark that the exceptional curves $L_{\alpha,\beta,\gamma}$
($\alpha$, $\beta$, $\gamma\in \mu_4$) are defined over $F^u$.
The quotient $G':=G/\langle u\rangle$ is isomorphic to the
dihedral group ${\mathfrak D}_4$;
generators $g:=\iota_a^3\iota_c$ and
$h:=\iota_b\iota_c^3\sigma$ satisfy $g^4=h^2=ghgh=e$.
We use the resolution of Proposition \ref{prop:efficient} to
identify classes in $H^1(G',M^u)$ with pairs
$(v,v')\in (M^u)^2$ satisfying
\begin{equation}
\label{eqn:cocyD4}
N_gv=N_hv'=0\qquad\text{and}\qquad N_{gh}v=N_{gh}v',
\end{equation}
modulo those of the form $(\Delta_gv,\Delta_hv)$.
Now a generator of $H^1(G',M^u)$ is the class of $(v_1,0)$
where
\begin{equation}
\label{eqn:v1}
v_1=(-1,0,1,0,0,0,0,0)=[L_{1,i,1}]-[L_{i,-1,-1}].
\end{equation}
Another representative for the same cohomology class is $(v_2,0)$ where
\begin{equation}
\label{eqn:v2}
v_2=(-1,0,-1,0,-1,-1,-2,2)=[L_{1,-1,i}]-[L_{i,i,i}].
\end{equation}

To produce an Azumaya algebra from one of these cocycles $(v_i,0)$
we must find rational equivalences that reflect the identities
\eqref{eqn:cocyD4}.
In fact, for each of the cycle representatives given in
\eqref{eqn:v1} and \eqref{eqn:v2}, the result of applying $N_{gh}$ is
equal to zero as a cycle.
So it remains only to find rational functions whose divisors are
$N_g$ applied to these cycle representatives.
For \eqref{eqn:v1}, a function that vanishes on
$L_{1,i,1}\cup L_{i,-i,-i}\cup L_{1,-i,1}\cup L_{i,i,-i}$
and has a simple pole along
$L_{i,-1,-1}\cup L_{1,-1,-i}\cup L_{i,1,-1}\cup L_{1,1,-i}$
is
$$f_1:=\frac{p(1+i)xz+iy^2-(1/2)w}{p(-1+i)xz+iy^2+(1/2)w}.$$
The corresponding rational equivalence for \eqref{eqn:v2} is
$$f_2:=\frac{p(1-i)xz+iy^2+(1/2)w}{p(1+i)xz-iy^2+(1/2)w}.$$
For $i=1$ and $2$ we have $f_i\,h(f_i)=1$, and the cocycle
$$(f_i,1,1)\in (F^u(S)^*)^3$$
determines an Azumaya algebra ${\mathfrak A}_i$ on $S$.

We claim ${\mathfrak A}_1$ and ${\mathfrak A}_2$ are equal in $\Br(S)$
and are:
\begin{itemize}
\item unramified at all points of $S(\Q_2)$;
\item ramified at all points of $S(\Q_{17})$;
\item unramified at all points of $S(\R)$.
\end{itemize}
The last of these claims is clear, since
$(f_i,1,1)\in {\mathbb S}^1\times\{1\}\times\{1\}$ at any point of $S(\R)$
(where ${\mathbb S}^1\subset\C^*$ denotes the unit circle) and this is a
connected subgroup of the group of cocycles, hence trivial in cohomology.

For the claim regarding $2$-adic points, we pause to discuss the
cohomology group $H^2({\mathfrak D}_4,\Q_2(i,\sqrt[4]{17})^*)$,
where generators act by
$$g\colon\begin{cases}i\mapsto i\\ \sqrt[4]{17}\mapsto i\sqrt[4]{17}\end{cases}
\qquad
h\colon\begin{cases}i\mapsto -i\\ \sqrt[4]{17}\mapsto i\sqrt[4]{17}\end{cases}
$$
Consider the diagram of field extensions, where labels indicate
fixed fields.
$$
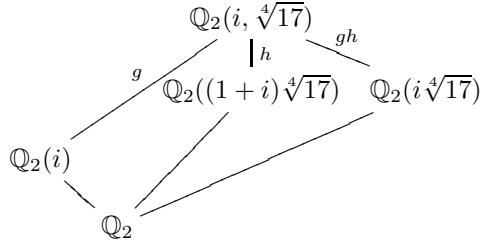
\xymatrix@R=10pt@C=5pt{
&&{\Q_2(i,\sqrt[4]{17})}\ar@{-}[ddll]_g\ar@{-}[d]^h\ar@{-}[dr]^{gh}\\
&&{\Q_2((1+i)\sqrt[4]{17})}\ar@{-}[ddl]&{\Q_2(i\sqrt[4]{17})}\ar@{-}[ddll]\\
{\Q_2(i)}\ar@{-}[dr]\\
&{\Q_2}
}$$
Now by the resolution for ${\mathfrak D}_4$ of Section \ref{sect:gr},
a $2$-cocycle is $(r,s,t)$ with
$$r\in \Q_2(i)^*,\qquad
s\in \Q_2((1+i)\sqrt[4]{17})^*, \qquad
t\in \Q_2(i\sqrt[4]{17})^*$$
satisfying
$Nr = Ns\,Nt$,
where in each instance, $N$ denotes the norm from the respective field
to $\Q_2$.
Coboundaries are triples
$$(N_g c , N_h d , N_{gh}(c/d))$$
for $c$, $d\in \Q_2(i,\sqrt[4]{17})^*$.

At every $2$-adic point of $S$, at least one of $f_1$ and $f_2$
is defined and takes one of the following values mod $32$:
\begin{equation}
\label{eqn:32}
\arraycolsep 4pt
\begin{array}{llllllll}
   1{+}0i&1{+}8i &1{+}16i &1{+}24i&25{+}4i&25{+}12i& 25{+}20i& 25{+}28i\\
   0{+}31i& 8{+}31i&16{+}31i&24{+}31i& 4{+}7i &12{+}7i&20{+}7i&28{+}7i\\
   31{+}0i& 31{+}24i& 31{+}16i&31{+}8i&7{+}28i&7{+}20i&7{+}12i&7{+}4i\\
   0{+}i &24{+}i &16{+}i  &8{+}i  &28{+}25i& 20{+}25i& 12{+}25i& 4{+}25i\\
\end{array}
\end{equation}
We claim that for any cocycle $(f,1,1)$ with
$f$ (necessarily in $\Z_2[i]$) taking one of the values mod $32$ listed
in \eqref{eqn:32}, there exists $c\in \Q_2(i,\sqrt[4]{17})^*$ with
$N_{gh} c = 1$ and $N_g c = f$, so
in particular, $(f,1,1)$ is a coboundary.
Indeed, the image of $N_g$ among $c\in \Z_2[i,\sqrt[4]{17}]$
satisfying $N_{gh} c=1$ is the set of
$f\in \Q(i)^*$ with $Nf=1$ and $f$ mod $32$ equal to some
value in the first row of \eqref{eqn:32}.
Also, there exists $c\in \Q_2(i,\sqrt[4]{17})^*$ with
$N_{gh} c = 1$ and $N_g c = i$.
Since norms are multiplicative, the claim follows.

The equality of ${\mathfrak A}_1$ and ${\mathfrak A}_2$ in $\Br(S)$
follows from having a function 
$r\in F^u(S)^*$, whose norm by $g$ is $f_2/f_1$ and whose norm by $gh$ is $1$.
Recall that $(v_2,0)$ equals $(v_1,0)$ in cohomology; explicitly this is
by $v_2-v_1=\Delta_g([L_{1,1,1}]+[L_{1,i,-i}])$.
Now $r$ can be taken to be a rational function vanishing on
$L_{1,-1,i}\cup L_{i,-1,-1}\cup g(L_{1,1,1})\cup g(L_{1,i,-i})$ with
a simple pole on
$L_{i,i,i}\cup L_{1,i,1}\cup L_{1,1,1}\cup L_{1,i,-i}$, scaled appropriately.

The $17$-adic analysis is simpler because
$\Q_{17}$ has $\sqrt{-1}$, and hence
we are reduced to analyzing norms for
$\Q_{17}\to\Q_{17}(\sqrt[4]{17})$.
Norms for this extensions are just
powers of $17$ times $4$-th powers in $\Z_{17}^*$.
Evaluating $f_1$ at points of $S(\Q_{17})$ and substituting
$\sqrt{-1}$ for $i$ yields the classes $8$ and $15$ mod $17$,
and these are not quartic residues.
\end{exam}

\begin{rema}
\label{rema:automate}
The analysis we have carried out in the examples could, in principle,
be carried out algorithmically in any of the arithmetic classes of
surfaces $S$.
We have verified that, except in two uninteresting cases (in which one
of $A$, $B$, and $C$ has to be a square),
the $2$-torsion subgroup of $\Br(S)/\Br(\Q)$ is generated by
the groups $H^1(\Z/2,M^H)$ as $H$ ranges over the index $2$ subgroups
of Galois group.
In the cases with $4$-torsion in $\Br(S)/\Br(\Q)$,
the analysis can proceed as in Example \ref{exam:obstr5}.
\end{rema}

\begin{rema}
\label{rema:end}
In Examples \ref{exam:obstr1} through \ref{exam:obstr5},
the surface $S$ always satisfies $\Pic(S)=\Z$.
So, in considering the cases of
del Pezzo surfaces of degree $2$ covered by previous results,
described in the introduction,
we are consistently avoiding the non-minimal surfaces of case (i).
Every surface that we are considering is
in some obvious ways, a double cover of
Ch\^atelet surfaces
(one can pass to invariants for any projective linear transformation
of $x$, $y$, and $z$
which is an involution preserving $Ax^4+By^4+Cz^4$).
But in every example,
the resulting Ch\^atelet surfaces satisfy the Hasse principle
(this can be seen by directly exhibiting rational points, combined with
appeal to \cite[Theorem B]{chateletsurfaces}).
In at least one case, Example~\ref{exam:obstr5},
it is easy to exclude the surface from 
being birational to a conic bundle, since
$\Br(S)/\Br(\Q)$, a birational invariant, is $2$-torsion for conic bundles.
\end{rema}

\section{Appendix: Cyclic Azumaya algebras on diagonal cubics}
In \cite{cubdiag}, there is
an analysis of the Brauer-Manin obstruction on
a diagonal cubic surface $S$, given by
\begin{equation}
\label{eqn:cubic}
Ax^3+By^3+Cz^3+Dt^3=0,
\end{equation}
with $A$, $B$, $C$, and $D$ positive integers.
Let $\theta=e^{2 \pi i/3}$; first of all,
$S(\Q)=\emptyset$ if and only if
$S(\Q(\theta))=\emptyset$, and hence it suffices to work over
the field $k:=\Q(\theta)$.
The analysis proceeds by constructing
Azumaya algebras that are split by a bicyclic extension of $k$
and computing local invariants.

Here we simplify the algorithm proposed in \cite{cubdiag} by constructing
{\em cyclic} Azumaya algebras on $S_k$ which generate $\Br(S_k)/\Br(k)$.
We use descent to exhibit the necessary cycles, as in
Example \ref{exam:obstr4}.

We start by making the following assumption:
\begin{gather}
\label{eqn:genericcubic}
\begin{split}
\sqrt[3]{A/B}\notin\Q,\ \sqrt[3]{A/C}\notin\Q,\ \ldots,\ \sqrt[3]{C/D}\notin\Q\\
\sqrt[3]{AB/CD}\notin\Q,\ \sqrt[3]{AC/BD}\notin\Q,\ \sqrt[3]{AD/BC}\notin\Q
\end{split}
\end{gather}
(in all other cases, the Hasse principle is known to hold).
Then we define
\begin{gather*}
\alpha=\sqrt[3]{B/A}\qquad \beta=\sqrt[3]{D/C}\qquad \gamma=\sqrt[3]{AD/BC}=
\alpha^{-1}\beta \\
\alpha'=\sqrt[3]{C/A}\qquad \beta'=\sqrt[3]{D/B}
\end{gather*}
We assume, further, that
$S(\Q_p)\ne\emptyset$ for all primes $p$.
Set $K=k(\gamma,\alpha)$; the assumption \eqref{eqn:genericcubic} implies
\begin{equation}
\label{eqn:nine}
[K:k]=9.
\end{equation}

We need notation for the following divisors on $S_{\bar k}$:
$$
\arraycolsep 0pt
\begin{array}{rclrclrcl}
L(i)&:&\begin{cases}
x{+}\theta^i\alpha y=0\\
z{+}\theta^i\beta t=0
\end{cases}&
L'(i)&:&\begin{cases}
x{+}\theta^i\alpha y=0\\
z{+}\theta^{i{+}1}\beta t=0
\end{cases}&
L''(i)&:&\begin{cases}
x{+}\theta^i\alpha y=0\\
z{+}\theta^{i{+}2}\beta t=0
\end{cases}
\end{array}
$$
and
$$
\arraycolsep 0pt
\begin{array}{rclrclrcl}
M(i)&:&\begin{cases}
x{+}\theta^i\alpha' z=0\\
y{+}\theta^{i{+}1}\beta' t=0
\end{cases}
\end{array}
$$
Define
$$L=L(0)+L(1)+L(2)\qquad\text{and}\qquad M=M(0)+M(1)+M(2).$$
Now $L+M$ is comprised of 6 pairwise disjoint lines; blowing these down
we have $S_{\bar k}\to \PP^2_{\bar k}$.
Take $\ell$ to be the class of a general line in $\PP^2_{\bar k}$, so
$$3\ell = -K_S+L+M.$$
By results in \cite{cubdiag}, we have
$$\Z/3=H^1(\Z/3,\Pic(S_{k(\gamma)}))\eqto \Br(S_k)/\Br(k),$$
generated by the class in $H^1(\Z/3,\Pic(S_{k(\gamma)}))$ of
$\ell-L$ or $\ell-M$ (where we use \eqref{eqn:cohocyclic} to
identify elements with cohomology classes).
In \cite{cubdiag}, the following procedure is proposed to
obtain a nontrivial Azumaya algebra on $S_k$:
\begin{itemize}
\item[(i)] Find a divisor $D$ defined over $k(\gamma)$ in the class
$\ell-L$ or $\ell-M$,
\item[(ii)] Find a function in $k(S)$ whose divisor is the union of $D$
and its Galois conjugates.
\end{itemize}
Unfortunately, the
classes in $\Pic(S_{k(\gamma)})$ of sums of lines defined over
$S_{k(\gamma)}$ fail to represent any nonzero elements of
$H^1(\Z/3,\Pic(S_{k(\gamma)}))$, and the further field
extension required to find suitable sums of lines accounts for much
of the complication of the analysis of \cite{cubdiag}.

We show that (i) can be carried out by solving a norm equation.
Then (ii) reduces to some linear algebra.
For (i), we start with the further field extension $k(\gamma)\to K$
and the divisor
$D:=L'(2)-L''(0)$
in class $\ell-M$ (cf.\ \cite{cubdiag}).
Denote by $\sigma$ the element of $\Gal(K/k(\gamma))$ which sends
$\alpha$ to $\theta\alpha$.
For the line bundle $\cO_{S_K}(D)$ to descend to $k(\gamma)$
we must supply
an isomorphism
$$\cO_{S_K}(L'(2)-L''(0))\stackrel{\xi}\lra \cO_{S_K}(L'(0)-L''(1))$$
satisfying
\begin{equation}
\label{eqn:cubiccocycle}
\sigma^2(\xi)\circ \sigma(\xi)\circ \xi = 1.
\end{equation}
Looking at the defining equations, we see $\xi$ must be of the form
$$\xi=\varepsilon \frac{z+\beta t}{x+\alpha y}$$
for some $\varepsilon\in k(\gamma)$.
Now the condition \eqref{eqn:cubiccocycle} is equivalent to
\begin{equation}
\label{eqn:cubicnorm}
N_{K/k(\gamma)}(\varepsilon)=-C/A.
\end{equation}
Concretely, if
$$\varepsilon=\lambda+\mu\alpha+\nu\alpha^2$$
with $\lambda$, $\mu$, $\nu\in k(\gamma)$, then \eqref{eqn:cubicnorm}
expands as
\begin{equation}
\label{eqn:cubicnormexpanded}
\lambda^3+\frac{B}{A}\mu^3+\frac{B^2}{A^2}\nu^3-3\frac{B}{A}\lambda\mu\nu=
-\frac{C}{A}.
\end{equation}
Equation \eqref{eqn:cubicnormexpanded}
has a solution, by the Hasse principle.
There is also an {\em a priori} bound on the size of some
solution \cite{siegel}.
An effective algorithm exists; see for example \cite{fjp}.
Algorithms from \cite{cohen} and \cite{fiekerthesis} have been
implemented in {\tt magma}.

Define $k'=k(\gamma)$.
By descent we have a line bundle $\cE$ on $S_{k'}$.
Also by descent, a rational section of $\cE$ is given by
\begin{align*}
f &= 1+\sigma^2(\xi)+\sigma(\xi)\sigma^2(\xi) \\
 &= \frac
{(x{+}\theta\alpha y)(x{+}\theta^2\alpha y)
+ \sigma^2\varepsilon(x{+}\theta\alpha y)(z{+}\theta^2\beta t)
+ \sigma\varepsilon\sigma^2\varepsilon(z{+}\theta\beta t)(z{+}\theta^2\beta t)}
{(x{+}\theta\alpha y)(x{+}\theta^2\alpha y)}.
\end{align*}
Then, with respect to local trivializations of $\cE$, the section $f$
has a simple pole on $L''(0)+L''(1)+L''(2)$ and vanishes to order one
along some cubic curve $C$.
Hence
$$C=-2L-M+4\ell$$
in $\Pic(S_{k'})$, and $C+K_S=-L+\ell$ is a divisor as desired.

We compute $C^2=1$ and $C\cdot K_S=-3$, which implies that its genus is 
zero, so $C$ is geometrically a twisted cubic.
Denoting by $g$ the numerator of $f$,
explicit defining equations of $C\subset S$ over $K$ are
$g=\sigma(g)=\sigma^2(g)=0$.
It is possible to express
$$g=g_0+g_1\alpha+g_2\alpha^2$$
for $g_0$, $g_1$, $g_2\in k'[x,y,z,t]$, and
after a bit of algebra we find
\begin{align*}
g_0 &= x^2 +\lambda x z + (B/A)\nu x t \gamma + \theta^2 (B/A) \mu y t \gamma
+\theta^2 (B/A) \nu y z \\
&\quad{}+ [\lambda^2 -(B/A)\mu\nu]z^2
+(B/A)(\lambda\nu-\mu^2)z t\gamma
+ (B/A) [ (B/A)\nu^2-\lambda\mu ]  t^2\gamma^2 \\
g_1 &= -xy +\theta^2 \mu x z + \theta^2 \lambda x t \gamma
+ \theta\lambda y z + \theta (B/A) \nu y t\gamma
+ [(B/A) \nu^2 - \lambda\mu] z^2 \\
&\quad{}+[(B/A) \mu\nu -\lambda^2]z t\gamma
+(B/A) (\mu^2-\lambda\nu) t^2\gamma^2 \\
g_2 &= \theta \nu x z + \theta \mu x t \gamma + y^2 + \mu y z
+ \lambda y t \gamma \\
&\quad{} + (\mu^2-\lambda\nu) z^2
+ [\lambda \mu -(B/A)\nu^2] z t \gamma
 + [\lambda^2-(B/A)\mu\nu] t^2 \gamma^2 
\end{align*}
Now $C$ is defined over $k'$ as a subvariety of $S$ by the equations
\begin{equation}
\label{eqn:g0g1g2}
g_0=g_1=g_2=0.
\end{equation}
In fact, we have
\begin{multline*}
g_0(Ax-A\lambda z-B\nu\gamma t)
+g_1(-B\nu z - B\mu\gamma t)
+g_2(By-B\mu z-B\lambda\gamma t) \\
= Ax^3+By^3+Cz^3+Dt^3
\end{multline*}
so \eqref{eqn:g0g1g2} defines $C$ over $k'$ as a subvariety of $\PP^3$.
We have completed task (i).

For task (ii), we claim there exist linear polynomials
$\ell_0$, $\ell_1$, $\ell_2\in k'[x,y,z,t]$
such that the polynomial
\begin{equation}\label{eqn:gell}
h=g_0\ell_0 + g_1\ell_1 + g_2\ell_2
\end{equation}
is in $k[x,y,z,t]$ and is not proportional to
$(Ax^3+By^3+Cz^3+Dt^3)$.
Knowing this, a modern linear algebra solver can effectively produce such
$\ell_0$, $\ell_1$, and $\ell_2$.
Then the division algebra generated over $k(S)$
by noncommuting variables $r$ and $s$ subject to relations
$$r^3 = AD/BC,\qquad s^3 = h/x^3,\qquad sr=\theta rs,$$
is the restriction of an Azumaya algebra over $S_k$ generating
$\Br(S_k)/\Br(k)$.

To justify the claim, notice first that there exists a rational
function on $S_k$ whose divisor is $3H-C-\rho C-\rho^2C$,
where $H$ is a hyperplane section and $\rho$ is a generator of
$\Gal(k'/k)$.
Next, by a dimension computation, we have an isomorphism
$$H^0(\PP^3_k,\cO(3))/\langle Ax^3+By^3+Cz^3+Dt^3\rangle\to H^0(S,3H)$$
so this rational function must be of the form $h/\ell^3$
(assuming that $H$ is defined by the vanishing of the linear form $\ell$).
Finally, a syzygy computation shows that $h$ can be expressed
in the form \eqref{eqn:gell}.
Indeed, \eqref{eqn:g0g1g2} defines $C$ in $\PP^3$, so we know
$\ell^dh$ lies in the ideal $(g_0,g_1,g_2)$ of $k'[x,y,z,t]$, for some $d$.
Suppose $d\ge 1$ and
$$\ell^dh = \sum_{i=0}^2g_ir_i,$$
with $r_i\in k'[x,y,z,t]$ for $i=0$, $1$, $2$.
Now it suffices to show that
there exist $s_0$, $s_1$, $s_2\in k'[x,y,z,t]$ such that
$\sum_i g_is_i=0$, and $\ell$ divides $r_i-s_i$ for each $i$; then
we have $\ell^{d-1}h=\sum_i g_i(r_i-s_i)/\ell$ and we can proceed inductively.
In other words, it suffices to show that the map on Koszul complexes
for $(g_0, g_1, g_2)$,
induced by the quotient map $k'[x,y,z,t]\to k'[x,y,z,t]/(\ell)$,
gives rise to a surjection on the first homology modules.
It is enough to verify this over the algebraic closure, and
we are reduced
to the case of $(g_0, g_1, g_2)$ defining the twisted cubic,
for which it is a standard computation.

\end{document}